\documentclass[11pt]{article}

\usepackage{latexsym, amssymb, amsmath, theorem, comment}

\makeatletter \@addtoreset{equation}{section}

\theoremstyle{change}
\theoremheaderfont{\scshape}

\newtheorem{thm}{Theorem}[section]
\newtheorem{prop}[thm]{Proposition}
\newtheorem{lem}[thm]{Lemma}

\newtheorem{defin}[thm]{Definition}

\theorembodyfont{\rmfamily}
\newtheorem{conj}[thm]{Conjecture}
\newtheorem{rem}[thm]{Remark}

\def\qed{\hfill{$\blacksquare$}\medskip}

\def\bp{\noindent{\it Proof.}\ }

\def\7#1{{\mathbb #1}}

\newcommand{\rarr}{\rightarrow}

\begin{document}

\title{The Mathieu conjecture for $SU(2)$ reduced to an abelian conjecture}

\author{Michael M\"uger\footnote{Radboud University, Nijmegen, The Netherlands, {\tt
mueger@math.ru.nl}}\ \  and Lars Tuset\footnote{Oslo Metropolitan University, Oslo, Norway, {\tt larst@oslomet.no}}}
\date{\today}
\maketitle

\begin{abstract}
We reduce the Mathieu conjecture for $SU(2)$ to a conjecture about moments of Laurent polynomials in
two variables with single variable polynomial coefficients.
\end{abstract}

\section{Introduction}

O.\ Mathieu conjectured \cite{mathieu} that for complex valued regular (=finite type) functions $f,g$
on any connected compact Lie group with normalized Haar integral $\int$, the vanishing of $\int f^n$ for all 
positive integers $n$ implies $\int f^n g =0$ for all large enough $n$. He then proved that this
conjecture implies O.-H.\ Keller's notorious Jacobian conjecture.\footnote{For lack of a better
concise name, we will follow \cite{FPYZ} in calling integrals like $\int f^ng$ `moments', since they
generalize classical moments of the form $\int_a^b x^n f(x)dx$.}

Motivated by work of Dings and Koelink \cite{dings} we reduce the Mathieu conjecture for $SU(2)$ to a  
conjecture of a more abelian nature about moments of Laurent polynomials in two variables with single variable
polynomial coefficients. Generalized further in the natural fashion, this `$xz$-conjecture' says that
if $f(x,z)=\sum_{\bf m}c_{\bf m}(x) z^{\bf m}$ is a Laurent polynomial in several $z$-variables with
polynomial coefficients $c_{\bf m}$ in several $x$-variables satisfying $\int f^n=0$ for all positive integer
$n$, where $z$ is integrated over the torus and $x$ over a cube, then $\bf 0$ is not in the convex hull of the 
set of multi-indices $\bf m$ for which $c_{\bf m}\neq 0$.  

In the absence of $x$-variables, our conjecture reduces to a result proven by Duistermaat and van der Kallen
\cite{dvdk} as part of their proof of the Mathieu conjecture in the abelian case. On the other hand, the
$xz$-conjecture with one $x$ variable and no $z$'s is known to hold, see \cite{MT} and references therein. 
For the moment, the $xz$-conjecture remains open already for one $z$ and one $x$. Towards the end of this
paper we explain that the natural inductive approach to proving it in this case fails due to the topological
`worm problem'. We also include a trivial generalization of the approach of Dings and Koelink to any connected  
compact Lie group.


\section{The $SU(2)$-case}
In what follows, we equip $\7T=\{z\in\7C\ | \ |z|=1\}$ with the normalized Haar measure
$\frac{1}{2\pi i}\,\frac{dz}{z}$ and $[0,1]$ with Lebesgue measure. Products $\7T^k\times[0,1]^l$ will carry
the obvious product measure.

Define maps
\begin{eqnarray*}
\alpha: && SU(2)\rarr\7C^4,\quad\begin{pmatrix} a & c\\ b & d\end{pmatrix}\mapsto(a,b,c,d),\\
   \beta: &&\7T^2\times[0,1]\rarr\7C^4, \quad (z_1,z_2,x)\mapsto((1-x)z_2,xz_1,-z_1^{-1}, z_2^{-1}).
\end{eqnarray*}
A function $f: SU(2)\rarr\7C$ is called regular if $f=P\circ\alpha$ for some $P\in\7C[a,b,c,d]$. (That $P$ is
not uniquely determined will not be an issue.) Denoting the space of regular functions by $R$, the maps
$\Lambda: \7C[a,b,c,d]\rarr R, P\mapsto P\circ\alpha$ and
$\Pi: \7C[a,b,c,d]\rarr C(\7T^2\times[0,1],\7C), P\mapsto P\circ\beta$ clearly are ring homomorphisms. 

\begin{lem} For each $P\in\7C[a,b,c,d]$ we have
\begin{equation}
  \int_{SU(2)}P\circ\alpha=\int_{\7T^2\times[0,1]} P\circ\beta. \label{eq-main}
\end{equation}  
\end{lem}

\bp It suffices to prove this for monomials $P=a^{n_1}b^{n_2}c^{n_3}d^{n_4}$. 
For the l.h.s.\ of (\ref{eq-main}), like Dings and Koelink \cite{dings} we use a classical
integration formula for $SU(2)$ \cite[Ch.\ III, Sect.\ 6.1]{V}:
\[ \int_{SU(2)} f =\frac{1}{16\pi^2}\int_0^{2\pi}\int_0^\pi\int_{-2\pi}^{2\pi}
  f(F(\phi,\theta,\psi))\sin\theta \,d\psi\,d\theta\,d\phi, \]
where
\[ F(\phi,\theta,\psi)=\begin{pmatrix} \cos\frac{\theta}{2}\,e^{\frac{i(\phi+\psi)}{2}} &
    i\sin\frac{\theta}{2}\,e^{\frac{i(\phi-\psi)}{2}} \\
    i\sin\frac{\theta}{2}\,e^{\frac{i(\psi-\phi)}{2}} & \cos\frac{\theta}{2}\,e^{\frac{-i(\phi+\psi)}{2}}
    \end{pmatrix}. \]
For $f=P=a^{n_1}b^{n_2}c^{n_3}d^{n_4}$ this becomes
\begin{eqnarray*} \int_{SU(2)}P\circ\alpha &=&
  \frac{i^{n_2+n_3}}{16\pi^2}\int_0^\pi\cos^{n_1+n_4}\frac{\theta}{2}\,\cdot\,\sin^{n_2+n_3}\frac{\theta}{2}\,\sin\theta d\theta\\
   && \int_0^{2\pi} e^{i\frac{\phi}{2}(n_1-n_2+n_3-n_4)} d\phi\ \int_{-2\pi}^{2\pi} e^{i\frac{\psi}{2}(n_1+n_2-n_3-n_4)} \,d\psi.
\end{eqnarray*}
Now, $\int_{-2\pi}^{2\pi}e^{i\frac{\psi}{2}(n_1+n_2-n_3-n_4)}d\psi=4\pi\delta_{n_1+n_2-n_3-n_4,0}$. When this
is non-zero, we have $n_1-n_4=n_3-n_2$, so that $n_1-n_2+n_3-n_4$ is
even, so that the integration $\int_0^{2\pi}\cdots d\phi$ gives a factor $2\pi\delta_{n_1-n_2+n_3-n_4,0}$. The
combination of $n_1+n_2-n_3-n_4=0$ and $n_1-n_2+n_3-n_4=0$ is equivalent to $n_1=n_4\wedge n_2=n_3$. Thus
\[ \int_{SU(2)} P\circ\alpha=\frac{(-1)^{n_2}\delta_{n_1,n_4}\delta_{n_2,n_3}}{2}
   \int_0^\pi\cos^{2n_1}\frac{\theta}{2}\ \,\sin^{2n_2}\frac{\theta}{2}\ \sin\theta \, d\theta. \]
With $x=\sin^2\frac{\theta}{2}$ we have
$\frac{dx}{d\theta}=2\sin\frac{\theta}{2}\cdot \cos\frac{\theta}{2}\cdot \frac{1}{2}=\frac{\sin\theta}{2}$,
thus $\sin\theta\,d\theta=2dx$, and with $\cos^2\frac{\theta}{2}=1-\sin^2\frac{\theta}{2}=1-x$ we have
\[ \int_{SU(2)} P\circ\alpha=(-1)^{n_2}\delta_{n_1,n_4}\delta_{n_2,n_3} \int_0^1 (1-x)^{n_1} x^{n_2} dx, \]
(where the $x$-integral is Euler's function $\beta(n_1+1,n_2+1)$). On the other hand,
\begin{eqnarray*} \int_{\7T^2\times[0,1]} P\circ\beta &=&
   \int_{\7T^2\times[0,1]} ((1-x)z_2)^{n_1}(xz_1)^{n_2}(-z_1^{-1})^{n_3}(z_2^{-1})^{n_4} \\
  &=& (-1)^{n_3} \int_0^1 (1-x)^{n_1}x^{n_2}dx\ \int_{\7T^2} z_1^{n_2-n_3}z_2^{n_1-n_4} \\
  &=& (-1)^{n_3}\delta_{n_1,n_4}\delta_{n_2,n_3}\int_0^1 (1-x)^{n_1} x^{n_2} dx,
\end{eqnarray*}
and comparing the two integrals completes the proof.
\qed

\begin{defin} Let $k,l\in\7N_0=\{0,1,\ldots\}$ and
$f\in\7C[x_1,\ldots,x_l, z_1,z_1^{-1},\ldots,z_k,z_k^{-1}]$. Considering $f$ as a Laurent  
polynomial $\sum_{{\bf m}\in\7Z^k} c_{\bf m}z^{\bf m}$ in $z_1,\ldots,z_k$ with coefficients
$c_{\bf m}\in\7C[x_1,\ldots,x_l]$, we define the spectrum of $f$ as
\[ \mathrm{Sp}(f)=\{ {\bf m}\in\7Z^k\ | \ c_{\bf m}\ne 0\}. \]
\end{defin}

\begin{conj} [$xz$-conjecture] Let $k,l\in\7N_0$ and
$f\in \7C[x_1,\ldots,x_l, z_1,z_1^{-1},\ldots,z_k,z_k^{-1}]$. If $\int_{[0,1]^l\times\7T^k} f^n=0$ for all
$n\in\7N$, then ${\bf 0}$ is not in the convex hull of $\mathrm{Sp}(f)\subset\7R^k$.
\end{conj}

\begin{rem} The conjecture is trivially true for $k=l=0$. For $l=0$, in which case the $c_{\bf m}$ are just
numbers, it was proven in \cite{dvdk}. For $k=0$ and all $l$ it was proven in \cite{FPYZ} and again by the
authors \cite{MT}, using ideas from the proof in \cite{dvdk} for $l=0, k=1$. To the best of the authors'
knowledge it is open for all other $(k,l)$. See Remark \ref{rem-worm} for comments on a failed attempt at
proving it for $k=l=1$.  
\end{rem}

\begin{thm} \label{theorem}
The case $k=2, l=1$  of the $xz$-conjecture implies the Mathieu conjecture for $SU(2)$.
\end{thm}

\bp Let $f$ be a regular function on $SU(2)$ such that $\int f^n=0$ for all $n\in\7N$. Pick $P\in\7C[a,b,c,d]$
such that $f=P\circ\alpha$. With the lemma, we have
\[ \int_{\7T^2\times[0,1]}(P\circ\beta)^n=\int_{\7T^2\times[0,1]}P^n\circ\beta
  =\int P^n\circ\alpha = \int (P\circ\alpha)^n =\int f^n=0 \quad\forall n\in\7N. \]
Since $P\circ\beta\in \7C[x,z_1,z_1^{-1},z_2,z_2^{-1}]$, the $xz$-conjecture (for $k=2, l=1$) implies that
${\bf 0}$ is not in the convex hull of $\mathrm{Sp}(P\circ\beta)\subset\7Z^2$. By a classical result, there is
a straight line in $\7R^2$ separating ${\bf 0}$ from $\mathrm{Sp}(P\circ\beta)$. This is equivalent to the
existence of ${\bf v}\in\7R^2\backslash\{{\bf 0}\}$ such that ${\bf v}\cdot {\bf m}\ge 1$ for all
${\bf m}\in\mathrm{Sp}(P\circ\beta)$. This implies ${\bf v}\cdot {\bf m}\ge n$ for all
${\bf m}\in\mathrm{Sp}((P\circ\beta)^n)$. As a consequence, $\mathrm{Sp}((P\circ\beta)^n)$ moves off to
infinity as $n\rarr\infty$, in the sense of becoming disjoint from  
every finite subset of $\7Z^2$ for $n$ large enough.

If now $g$ is another regular function on $SU(2)$ and $Q$ a polynomial such that $g=Q\circ\alpha$, for each
$n\in\7N$ we have
\[ \int f^ng=\int(P\circ\alpha)^n(Q\circ\alpha)=
\int (P^nQ)\circ\alpha=\int_{\7T^2\times[0,1]}(P^nQ)\circ\beta
    =\int_{\7T^2\times[0,1]}(P\circ\beta)^n(Q\circ\beta). \]
By the above, $\mathrm{Sp}((P\circ \beta)^n)$ will be disjoint from the finite set $-\mathrm{Sp}(Q\circ\beta)$
for all large enough $n$, so that the $z$-integrations give zero for all $x\in[0,1]$. Thus
$\int f^ng=0$, proving the Mathieu conjecture for $SU(2)$.
\qed

\begin{rem} \label{rem-worm}
We briefly report on a failed attempt to prove the $xz$-conjecture for one $x$ and one $z$ by adapting the
approach to Laurent polynomials in one variable $z$ pursued by Duistermaat and van der Kallen. 
Namely for $f=f(z;x)$ with $0$ in the convex hull of $\mathrm{Sp}(f)$, consider the generating function
$$
F(t)=\sum_{n=1}^\infty t^{n-1}\int_0^1 dx\int_{\mathbb{T}} f^n
=(2\pi i)^{-1}\int_0^1 dx\int_C\frac{f(z;x)}{1-tf(z;x)}\frac{dz}{z} , 
$$
which defines a holomorphic function for $|t|$ small. In analogy to \cite{dvdk}, as $\lim_{z\to 0}
f(z)=\infty$, the residue theorem and L'Hopital's rule tell us that for such $t$, we have 
$$
F(t)=-\frac{1}{t}-\sum_j\frac{1}{t^2}
\int_0^1 \frac{dx}{f'(\zeta_j ;x)\zeta_j} ,
$$
where $\zeta_j (z)$ are the $x$-dependent solutions of $f(\zeta_j ; x)=\tau$ with $|\zeta_j |<1$.  Of course this requires that $\tau =1/t$ is not a critical value of $f$, so that we don't get zero in the denominators. One would like then to analytically extend these $\zeta_j$'s along some curve $C'$ that avoids the critical values of $f$, and arguing that the functions under the latter integrals are bounded, conclude that the corresponding $F$ extends analytically to a function that is not identically zero by looking at its behavior 
as $\tau\to 0$. In the final round one needs then to carefully discuss the contributions to the integrals of the various cases $\zeta_j (\tau)$ as $\tau\to 0$ having in mind that $\zeta_j (\tau)$ might well converge to critical points of $f$. The hope then is that $-1/t$ will be the dominating term in $F$ in the limit, whence $F\neq 0$, showing that the moments of $f$ cannot all vanish.

But here the problem arises  that we don't have $C'$ even for simple $f$'s. Indeed, the critical values of 
$f(z;x)=c_{-1} (x)z^{-1} +c_0 (x)+c_1 (x)z$
are given by $\tau_{\pm} (x) =c_0 (x)\pm 2\sqrt {c_{-1}(x)c_1 (x)}$, and the specific choice of polynomials $c_j (x)$ given by
$c_1 (x)=c_{-1}(x)=2x-1+i(1-(2x-1)^2 )$ and $c_0 (x)=2x-1-i(1-(2x-1)^2 )$
produces curves or `worms' $\tau_{\pm}([0,1])$ that enclose the origin completely, thus preventing any curve
to reach the origin from infinity.    
\end{rem}


\section{The Dings-Koelink approach}

Let $G$ be a connected compact Lie group with maximal torus $T$, say of dimension $r$. Let $\hat{G}$ be the
set of equivalence classes of irreducible unitary
representations of $G$, and let $V_\sigma$ be the $G$-module 
of $\sigma\in\hat{G}$. Decompose 
$V_\sigma =\oplus_{m\in\hat{T}}V_{\sigma ,m}$
as a module over $T$. Let $f$ be a regular function on $G$. Then by definition we may write
$$
f=\sum_{\sigma\in\hat{G}} \mathrm{Tr}_{V_{\sigma}} (A_{\sigma}\pi_{\sigma}(\cdot))
=\sum_{\sigma\in\hat{G}} \sum_{m,m'\in\hat{T}}\mathrm{Tr}_{V_{\sigma}} (A_{\sigma , m,m'}\pi_{\sigma ,m,m'}(\cdot))
$$
for only finitely many non-zero complex quadratic matrices $A_{\sigma}$ each of 
size $\dim (V_\sigma )$. 
Consider the `spectrum' of $f$ to be 
$$
X_f=\{ (m,m')\in\hat{T}\times\hat{T}\ |\ A_{\sigma ,m,m'}\neq 0\}\subset\mathbb{R}^r\times\mathbb{R}^r .
$$ 
Dings and Koelink then showed by using the multinomial formula and left- and right actions of $T$ on $G$, that if $(0,0)$ is not in the convex hull of $X_f$, then
the moments $\int_G f^n (s)\, ds$ vanish for all non-negative integers $n$, where $ds$ is the Haar measure on $G$. By using the same trick once more, they also showed that if the converse (their Conjecture 4.1) of the previous statement holds
for $G$, then the Mathieu conjecture holds for $G$.
Thus it remains to show their conjecture:
\medskip

{\it If all the moments of $f$ vanish, then $(0,0)$ is not in the 
convex hull of $X_f$.}

\medskip
   
One might of course ask what the relation between $X_f$ and $\mathrm{Sp}(f)$ is when $G=SU(2)$, but we won't discuss that here.

\vspace{.5cm}
 
\noindent{\it Acknowledgment.} The authors thank the referees for comments simplifying the exposition. 
M.\ M.\ thanks L.\ T.\ and OsloMet for hospitality during two
months spent there in 2019. His expenses were shared by OsloMet and Radboud University.

\end{document}